# Lyapunov Characterizations of Input to Output Stability

Eduardo Sontag*
Dept. of Mathematics
Rutgers University
New Brunswick, NJ 08903
sontag@control.rutgers.edu

Yuan Wang†
Dept. of Mathematics
Florida Atlantic University
Boca Raton, FL 33431
ywang@math.fau.edu

**Abstract**

This paper presents necessary and sufficient characterizations of several notions of input to output stability. Similar Lyapunov characterizations have been found to play a key role in the analysis of the input to state stability property, and the results given here extend their validity to the case when the output, but not necessarily the entire internal state, is being regulated.

*Running Head.* Lyapunov Characterizations of Input to Output Stability

## 1 Introduction

This paper concerns itself with systems with outputs of the general form

$$\dot{x}(t) = f(x(t), u(t)), \quad y(t) = h(x(t)), \tag{1}$$

where $f : \mathbb{R}^n \times \mathbb{R}^m \to \mathbb{R}^n$ and $h : \mathbb{R}^n \to \mathbb{R}^p$ are both locally Lipschitz continuous, $f(0,0) = 0$, and $h(0) = 0$. In the paper [27] (see also [26]), the authors introduced several notions of output stability for such systems. All these notions serve to formalize the idea of a "stable" dependence of outputs $y$ upon inputs (which may be thought of as disturbances, actuator or measurement errors, or regulation signals). They differ in the precise formulation of the decay estimates and the overshoot, or transient behavior, characteristics of the output. Among all of them, the one of most interest is probably the one singled out for the name *input to output stability*, or IOS, for short.

Our main theorem in this paper provides a necessary and sufficient characterization of the IOS property in terms of Lyapunov functions. In the process of obtaining this characterization, we derive as well corresponding results for the variants of IOS discussed in [27]. (The relationships between those variants, shown in [27], play a role in our proofs, but otherwise the two papers are independent of each other.)

In the very special case when $y = x$, our concepts all reduce to the input to state stability (ISS) property. Much of ISS control design, cf. [3, 4, 6, 7, 8, 9, 10, 14, 16, 17, 28], relies upon the Lyapunov characterizations first obtained in [13, 23]. Thus, it is reasonable to expect a similar impact from the results given here for the more general case.

In order to review the different i/o stability concepts, let us make the following notational conventions. Euclidean norms will be denoted as $|x|$, and $\|u\|$ denotes the $L_\infty^m$-norm (possibly

*This research was supported in part by US Air Force Grant F49620-98-1-0242
†This research was supported in part by NSF Grants DMS-9457826 and DMS-9403924

infinite) of an input $u$ (i.e., a measurable and locally essentially bounded function $u : \mathcal{I} \to \mathbb{R}^m$, where $\mathcal{I}$ is a subinterval of $\mathbb{R}$ which contains the origin; if we do not specify the domain $\mathcal{I}$ of an input $u$, we mean implicitely that $\mathcal{I} = \mathbb{R}_{\geq 0}$). For each initial state $\xi \in \mathbb{R}^n$ and input $u$, we let $x(\cdot, \xi, u)$ be the unique maximal solution of the initial value problem $\dot{x} = f(x, u)$, $x(0) = \xi$, and write the corresponding output function $h(x(t, \xi, u))$ simply as $y(\cdot, \xi, u)$. Given a system with control-value set $\mathbb{R}^m$, we often consider the same system but with controls restricted to take values in some subset $\Omega \subseteq \mathbb{R}^m$; we use $\mathcal{M}_\Omega$ for the set of all such controls. As usual, by a $\mathcal{K}$ function we mean a function $\gamma : [0, \infty) \to [0, \infty)$ that is strictly increasing and continuous, and satisfies $\gamma(0) = 0$, by a $\mathcal{K}_\infty$ function one that is in addition unbounded, and we let $\mathcal{KL}$ be the class of functions $[0, \infty)^2 \to [0, \infty)$ which are of class $\mathcal{K}$ on the first argument and decrease to zero on the second argument. When we state the various properties below, we always interpret the respective estimates as holding for all inputs $u$ and for all initial states $\xi \in \mathbb{R}^n$.

With these conventions, we say that a system (1) is *uniformly bounded input bounded state stable*, and write "UBIBS" for short, if it is forward complete (for every initial state $\xi$ and input $u$, the solution $x(t, \xi, u)$ is defined for all $t \geq 0$) and, for some function $\sigma$ of class $\mathcal{K}$, the following estimate holds for all solutions:

$$|x(t, \xi, u)| \leq \max\{\sigma(|\xi|), \sigma(\|u\|)\}, \qquad \forall t \geq 0. \tag{2}$$

As argued in [27], stability of internal signals is a routine constraint in regulation problems. The following four output stability properties were discussed in [27]. A UBIBS system is:

- IOS, or *input to output stable*, if there exist a $\mathcal{KL}$-function $\beta$ and a $\mathcal{K}$-function $\gamma$ such that

$$|y(t, \xi, u)| \leq \beta(|\xi|, t) + \gamma(\|u\|), \qquad \forall t \geq 0 \tag{3}$$

(the term $\gamma(\|u\|)$ can be replaced by the norm of the restriction to past inputs $\gamma(\|u\|_{[0,t]})$, and the sum could be replaced by a "max" or two analogous terms);

- OLIOS, or *output-Lagrange input to output stable* if it is IOS and, in addition, there exist some $\mathcal{K}$-functions $\sigma_1, \sigma_2$ such that

$$|y(t, \xi, u)| \leq \max\{\sigma_1(|h(\xi)|), \sigma_2(\|u\|)\}, \qquad \forall t \geq 0; \tag{4}$$

- SIIOS, or *state-independent input to output stable* if there exist some $\beta \in \mathcal{KL}$ and some $\gamma \in \mathcal{K}$ such that

$$|y(t, \xi, u)| \leq \beta(|h(\xi)|, t) + \gamma(\|u\|), \qquad \forall t \geq 0; \tag{5}$$

- ROS, or *robustly output stable* (ROS) if there are a smooth $\mathcal{K}_\infty$-function $\lambda$ and a $\beta \in \mathcal{KL}$ such that
$$|y_\lambda(t, \xi, d)| \leq \beta(|\xi|, t), \qquad \forall t \geq 0, \tag{6}$$
where we denote by $y_\lambda(\cdot, \xi, d)$ the output function of the system
$$\dot{x} = g(x, d) := f(x, d\lambda(|y|)), \quad y = h(x), \tag{7}$$
and we understand the estimate (6) as holding for all $d \in \mathcal{M}_\mathcal{B}$, where $\mathcal{B} = \{|\mu| \leq 1\} \subset \mathbb{R}^m$.

The last concept corresponds to the preservation of output stability under output feedback with "robustness margin" $\lambda$. It was shown in [27] that SIIOS $\Rightarrow$ OLIOS $\Rightarrow$ IOS $\Rightarrow$ ROS, and no converses hold. We now introduce the associated Lyapunov concepts.

**Definition 1.1** With respect to the system (1), a smooth function $V : \mathbb{R}^n \to \mathbb{R}_{\geq 0}$ is:

- an IOS-*Lyapunov function* if there exist $\alpha_1, \alpha_2 \in \mathcal{K}_\infty$ such that

$$\alpha_1(|h(\xi)|) \leq V(\xi) \leq \alpha_2(|\xi|), \qquad \forall \xi \in \mathbb{R}^n, \qquad (8)$$

  and there exist $\chi \in \mathcal{K}$ and $\alpha_3 \in \mathcal{KL}$ such that

$$V(\xi) \geq \chi(|\mu|) \;\Rightarrow\; DV(\xi)f(\xi,\mu) \leq -\alpha_3(V(\xi), |\xi|), \qquad \forall \xi \in \mathbb{R}^n, \forall \mu \in \mathbb{R}^m; \qquad (9)$$

- an OLIOS-*Lyapunov function* if there exist $\alpha_1, \alpha_2 \in \mathcal{K}_\infty$ such that

$$\alpha_1(|h(\xi)|) \leq V(\xi) \leq \alpha_2(|h(\xi)|), \qquad \forall \xi \in \mathbb{R}^n, \qquad (10)$$

  and there exist $\chi \in \mathcal{K}$ and $\alpha_3 \in \mathcal{KL}$ such that (9) holds;

- an SIIOS-*Lyapunov function* if there exist $\chi \in \mathcal{K}$ and $\alpha_3 \in \mathcal{K}$ such that

$$V(\xi) \geq \chi(\mu) \Rightarrow DV(\xi)f(\xi,\mu) \leq -\alpha_3(V(\xi)), \qquad \forall \xi \in \mathbb{R}^n, \forall \mu \in \mathbb{R}^m, \qquad (11)$$

  and there exist $\alpha_1, \alpha_2 \in \mathcal{K}_\infty$ such that (10) holds;

- an ROS-*Lyapunov function* if there exist $\chi \in \mathcal{K}$ and $\alpha_3 \in \mathcal{KL}$ such that

$$|h(\xi)| \geq \chi(|\mu|) \;\Rightarrow\; DV(\xi)f(\xi,\mu) \leq -\alpha_3(V(\xi), |\xi|), \qquad \forall \xi \in \mathbb{R}^n, \forall \mu \in \mathbb{R}^m, \qquad (12)$$

  and there exist $\alpha_1, \alpha_2 \in \mathcal{K}_\infty$ such that (8) holds. □

Observe that, if an estimate (8) holds, then (12) is implied by (9), in the sense that if $\chi$ and $\alpha_1$ are as in the former, then $\widetilde{\chi} := \alpha_1^{-1} \circ \chi$ can be used as "$\chi$" for the latter. Note also that, provided that (10) holds, condition (9) is equivalent to the existence of $\chi \in \mathcal{K}$ and $\alpha_3 \in \mathcal{KL}$ so that

$$|h(\xi)| \geq \chi(|\mu|) \;\Rightarrow\; DV(\xi)f(\xi,\mu) \leq -\alpha_3(V(\xi), |\xi|).$$

Our main results can be summarized as follows.

**Theorem 1** *A* UBIBS *system is:*

1. IOS *if and only if it admits an* IOS-*Lyapunov function;*
2. OLIOS *if and only if it admits an* OLIOS-*Lyapunov function;*
3. ROS *if and only if it admits an* ROS-*Lyapunov function; and*
4. SIIOS *if and only if it admits an* SIIOS-*Lyapunov function.*

The proofs are provided in Section 4.

## 2 Remarks on Rates of Decrease

In properties (9) and (12), the decay rate of $V(x(t))$ depends on the state and on the value of $V(x(t))$. The main role of $\alpha_3$ is to allow for slower convergence if $V(x(t))$ is very small or if $x(t)$ is very large. We first note two simplifications.

**Remark 2.1** Inequality (9) holds for some $\alpha_3 \in \mathcal{KL}$ if and only if there exist $\mathcal{K}$-functions $\kappa_1, \kappa_2$ such that

$$V(\xi) \geq \chi(|\mu|) \;\Rightarrow\; DV(\xi)f(\xi,\mu) \leq -\frac{\kappa_1(V(\xi))}{1+\kappa_2(|\xi|)} \tag{13}$$

for all $\xi \in \mathbb{R}^n$ and all $\mu \in \mathbb{R}^m$. This follows from Lemma A.2, proved in the Appendix. A similar remark applies to (12). □

**Remark 2.2** Suppose $V$ is an IOS-Lyapunov function for the system satisfying (8) with some $\alpha_1, \alpha_2 \in \mathcal{K}_\infty$ and (13) with some $\chi, \kappa_1, \kappa_2 \in \mathcal{K}$. By the proof of Lemma 11 together with Lemma 12 in [16], one sees that there exists a $C^1$ $\mathcal{K}_\infty$-function $\rho$ such that $\rho'(s)\kappa_1(s) \geq \rho(s)$ for all $s \geq 0$. Let $W = \rho \circ V$. Then $W$ is a $C^1$ function satisfying the following:

$$\rho(\alpha_1(|h(\xi)|)) \leq W(\xi) \leq \rho(\alpha_2(|\xi|)), \qquad \forall \xi \in \mathbb{R}^n,$$

and

$$W(\xi) \geq \chi_1(|\mu|) \;\Rightarrow\; DW(\xi)f(\xi,\mu) \leq -\frac{W(\xi)}{1+\kappa_2(|\xi|)} \tag{14}$$

for all $\xi \in \mathbb{R}^n$ and all $\mu \in \mathbb{R}^m$, where $\chi_1 = \rho \circ \chi \in \mathcal{K}$. This shows that if a system admits an IOS-Lyapunov function, then it admits one satisfying inequality (14). A similar remark applies to (12). □

Obviously, a function which satisfies a decay estimate of the stronger form

$$V(\xi) \geq \chi(|\mu|) \;\Rightarrow\; DV(\xi)f(\xi,\mu) \leq -\alpha(V(\xi)) \tag{15}$$

for some $\chi, \alpha \in \mathcal{K}$ is in particular an IOS Lyapunov function. It is thus natural to ask if there always exists, for an IOS system, a function with this stronger property. We now show, by means of an example, that such functions do not in general exist. Consider for that purpose the following two-dimensional single-input system:

$$\dot{x}_1 = 0, \quad \dot{x}_2 = -\frac{2x_2 + u}{1 + x_1^2}, \quad y = x_2. \tag{16}$$

This system is IOS, because with $V(x) := x_2^2$, it holds that

$$V(\xi) \geq \mu^2 \;\Rightarrow\; DV(\xi)f(\xi,\mu) = -2x_2 \frac{2x_2 + u}{1 + x_1^2} \leq -\frac{2V(\xi)}{1 + x_1^2}.$$

Namely, $V$ is an IOS-Lyapunov function for the system.

Suppose that system (16) would admit an IOS-Lyapunov function $W$ with a decay estimate as in (15), i.e., there exist some $\chi, \alpha \in \mathcal{K}$ such that

$$W(\xi) \geq \chi(|\mu|) \;\Rightarrow\; DW(\xi)f(\xi,\mu) \leq -\alpha(W(\xi)). \tag{17}$$

Without loss of generality, we may assume that $\chi \in \mathcal{K}_\infty$. In particular, we have that

$$DW(\xi)f(\xi, -\chi^{-1}(W(\xi))) \leq -\alpha(W(\xi)) \tag{18}$$

for all $\xi \in \mathbb{R}^2$. Fix any $\xi_1 \in \mathbb{R}$, and consider the one-dimensional differential equation

$$\dot{x}_2 = -\frac{2x_2 - \chi^{-1}(W(\xi_1, x_2))}{1 + \xi_1^2}. \tag{19}$$

Since $W(\xi_1, x_2(t)) \to 0$ (because of (18)) and as $\alpha_1(|\xi_2|) \leq W(\xi_1, \xi_2)$ for all $\xi$ (for some $\alpha_1 \in \mathcal{K}$), it follows that $x_2(t) \to 0$ as $t \to \infty$. This implies that $W(\xi_1, \xi_2) < 2\xi_2$ for all $\xi_1 \in \mathbb{R}$ and $\xi_2 > 0$. Together with (17), this implies that there exists some $\beta \in \mathcal{KL}$ such that, for every trajectory of (16) with $u(t) \equiv 0$, it holds that

$$|x_2(t)| \leq \beta(|x_2(0)|, t),$$

for all $\xi = (x_1(0), x_2(0))$ such that $x_2(0) > 0$. This is impossible, as it can be seen that, when $u(t) \equiv 0$, $x_2(t) = x_2(0)e^{-2t/(1+(x_1(0))^2)}$ whose decay rate depends on both $x_2(0)$ and $x_1(0)$.

Observe that, if we let $U(\xi_1, \xi_2) := [(1 + \xi_1^2)|\xi_2|]^{(1+\xi_1^2)}$, then one obtains the following estimate:

$$|\xi_2| \geq |\mu| \Rightarrow DU(\xi)f(\xi, \mu) \leq -U(\xi) \tag{20}$$

for all $\xi_1 \in \mathbb{R}$, $\xi_2 \neq 0$, and all $\mu \in \mathbb{R}$. (The function $U$ is not smooth on the set where $U(\xi) = 0$, but using a routine smoothing argument, one may easily modify $U$ to get a smooth Lyapunov function.) This $U$ is not an example of a $W$ as here (which, in any case, we know cannot exist), because (20) only means that $U$ is an ROS-Lyapunov function, not necessarily an IOS-Lyapunov function (since the comparison is between $|\xi_2|$ and $|\mu|$ rather than between a function of $U$ and $|\mu|$).

Finally, we observe that property (9) in the IOS-Lyapunov definition may be rephrased as follows:

$$V(\xi) > \widetilde{\chi}(|\mu|) \Rightarrow DV(\xi)f(\xi, \mu) < 0, \qquad \forall \xi \in \mathbb{R}^n, \forall \mu \in \mathbb{R}^m \tag{21}$$

where $\widetilde{\chi}(s) := \rho\,\chi(s)$ (for any arbitrary chosen $\rho \in (0, 1)$). This statement is obviously implied by (9). Conversely, if $V$ satisfies this property, then there is an $\alpha \in \mathcal{KL}$ so that (9) holds; this follows from Lemma A.5 given in the Appendix.

## 3 Uniform Stability Notions

There is a key technical result which underlies the proofs of all our converse Lyapunov theorems. It requires yet another set of definitions, which correspond to stability uniformly on all "disturbance" inputs.

**Definition 3.1** A system (1) is *uniformly output stable* with respect to inputs in $\mathcal{M}_\Omega$, where $\Omega$ is a compact subset of $\mathbb{R}^m$, if:

- it is forward complete, and
- there exists a $\mathcal{KL}$-function $\beta$ such that

$$|y(t, \xi, u)| \leq \beta(|\xi|, t), \qquad \forall t \geq 0 \tag{22}$$

holds for all $u$ and all $\xi \in \mathbb{R}^n$.

If, in addition, there exists $\sigma \in \mathcal{K}$ such that

$$|y(t, \xi, u)| \leq \sigma(|h(\xi)|), \qquad \forall\, t \geq 0 \tag{23}$$

holds for all trajectories of the system with $u \in \mathcal{M}_\Omega$, then the system is *output-Lagrange uniformly output stable* with respect to inputs in $\mathcal{M}_\Omega$. Finally, if strengthens (22) to

$$|y(t, \xi, u)| \leq \beta(|h(\xi)|, t), \qquad \forall\, t \geq 0 \tag{24}$$

holding for all trajectories of the system with $u \in \mathcal{M}_\Omega$, then the system is *state-independent uniformly output stable* with respect to inputs in $\mathcal{M}_\Omega$.

**Theorem 2** *Let $\Omega$ be a compact subset of $\mathbb{R}^m$, and suppose that a system (1) is uniformly output stable with respect to inputs in $\mathcal{M}_\Omega$. Then the system admits a smooth Lyapunov function $V$ satisfying the following properties:*

- *there exist $\alpha_1, \alpha_2 \in \mathcal{K}_\infty$ such that*

$$\alpha_1(|h(\xi)|) \leq V(\xi) \leq \alpha_2(|\xi|), \qquad \forall\, \xi \in \mathbb{R}^n, \tag{25}$$

- *there exists $\alpha_3 \in \mathcal{KL}$ such that*

$$DV(\xi) f(\xi, \mu) \leq -\alpha_3(V(\xi), |\xi|), \qquad \forall\, \xi \in \mathbb{R}^n,\ \forall\, \mu \in \Omega. \tag{26}$$

*Moreover, if the system is output-Lagrange uniformly output stable with respect to inputs in $\mathcal{M}_\Omega$, then (25) can be strengthened to*

$$\alpha_1(|h(\xi)|) \leq V(\xi) \leq \alpha_2(|h(\xi)|), \qquad \forall\, \xi, \tag{27}$$

*for some $\alpha_1, \alpha_2 \in \mathcal{K}$.*

*Finally, if the system is state-independent uniformly output stable with respect to inputs in $\mathcal{M}_\Omega$, then (25) can be strengthened to (27) and also (26) can be strengthened to:*

$$DV(\xi) f(\xi, \mu) \leq -\alpha_4(V(\xi)), \qquad \forall\, \xi \in \mathbb{R}^n,\ \forall\, \mu \in \Omega. \tag{28}$$

*for some $\alpha_4 \in \mathcal{K}$.*

The proof of this theorem will be postponed until Section 4.6.

# 4 Proof of Theorem 1

We prove the various parts of the theorem after establishing a preliminary simple fact.

## 4.1 A Small-Gain Lemma for Output-Lagrange Stability

**Lemma 4.1** *For every system which satisfies (4), there exist a $\mathcal{K}$-function $\sigma$ and a $\mathcal{K}_\infty$-function $\lambda$ such that, for any $\xi$ and any $u$, if $|u(t)| \leq \lambda(|y(t, \xi, u)|)$ for almost all $t \geq 0$, then*

$$|y(t, \xi, u)| \leq \sigma(|h(\xi)|), \qquad \forall\, t \geq 0. \tag{29}$$

*Proof.* Assume system (1) satisfies (4) for some $\sigma_1, \sigma_2 \in \mathcal{K}$. We may suppose, without loss of generality, that $\sigma_2 \in \mathcal{K}_\infty$ and $\sigma_1(s) \geq s$. Let

$$\lambda(s) = \sigma_2^{-1}(\sigma_1^{-1}(s/2)/2).$$

Below we show that (29) holds for such a choice of $\lambda$, and $\sigma = \sigma_1$. Pick any $\xi$ and $u$ such that $|u(t)| \leq \lambda(|y(t, \xi, u)|)$. Let $x(t)$ denote the corresponding trajectory, and $y(t) = h(x(t))$.

Consider the case $h(\xi) \neq 0$. It is enough to show that the supremum $t_1$ of the times $s \geq 0$ for which $\gamma(|u(t)|) \leq |h(\xi)|/2$ a.e. on $[0, s]$ is infinite, since this together with (4) will imply that $|y(t)| \leq \sigma_1(|h(\xi)|)$ for all $t \geq 0$. If it were the case that $t_1 < \infty$ then, since $|u(t)| \leq |h(\xi)|/2$ almost everywhere on $[0, t_1]$, it follows from (4) that $|y(t)| \leq \sigma_1(|h(\xi)|)$ for all $t \in [0, t_1]$. By continuity, there exists some $\delta > 0$ such that $|y(t)| \leq 2\sigma_1(|h(\xi)|)$ for all $t \in [0, t_1 + \delta]$. Consequently,

$$\sigma_2(|u(t)|) \leq \sigma_2(\lambda(|y(t)|)) \leq \sigma_2(\lambda(2\sigma_1(|h(\xi)|))) \leq |h(\xi)|/2$$

for almost all $t \in [0, t_1 + \delta]$. This contradicts the definition of $t_1$. Hence, $t_1 = \infty$.

The case when $h(\xi) = 0$ is similar. Suppose that $|u(t)| \leq \lambda(|y(t)|)$ almost everywhere, but there is some $t_1 > 0$ and some $\varepsilon > 0$ such that $|y(t_1)| \geq \varepsilon$. Then $t_2 = \inf\{t : |y(t)| \geq \sigma_1^{-1}(\varepsilon/2)\} \in (0, t_1)$, and $|y(t_2)| = \sigma_1^{-1}(\varepsilon/2)$. Using the first part of the proof with the initial time $t_1$, we obtain that $|y(t)| \leq \sigma_1(|y(t_2)|) \leq \varepsilon/2$ for all $t \geq t_2$, contradicting the assumption that $|y(t_1)| \geq \varepsilon$. It follows that $y(t) \equiv 0$. ■

## 4.2 Proof of Theorem 1, Part 1.

*Necessity.* Consider an OLIOS system (1). By Lemma 4.1, there exist $\lambda_1 \in \mathcal{K}_\infty$ and $\sigma \in \mathcal{K}$ such that, for any $\xi$ and any $u$ such that $|u(t)| \leq \lambda_1(|y(t, \xi, u)|)$ a.e., (29) holds. Since the system is OLIOS, and in particular, IOS, and as shown in [27], any IOS system is necessarily also ROS, there exists some smooth $\mathcal{K}_\infty$-function $\lambda_2$ such that the system

$$\dot{x} = f(x, d\lambda_2(|y|)), \quad y = h(x), \tag{30}$$

where $d \in \mathcal{M}_\mathcal{B}$, is forward complete, and there exists some $\beta \in \mathcal{KL}$ such that, for all trajectories $x_{\lambda_2}(t, \xi, u)$ with the output functions $y_{\lambda_2}(t, \xi, u)$, it holds that

$$\left|y_{\lambda_2}(t, \xi, d)\right| \leq \beta(|\xi|, t), \qquad \forall\, t \geq 0,\ \forall\, \xi \in \mathbb{R}^n,\ \forall\, d \in \mathcal{M}_\mathcal{B}.$$

Let $\lambda_3(s) = \min\{\lambda_1(s), \lambda_2(s)\}$, and let $\lambda(\cdot)$ be any smooth $\mathcal{K}_\infty$-function so that $\lambda(s) \leq \lambda_3(s)$ for all $s$. Then, for the system

$$\dot{x} = f(x, d\lambda(|y|)), \quad y = h(x), \tag{31}$$

where $d \in \mathcal{M}_\mathcal{B}$, it holds that

$$\left|y_\lambda(t, \xi, d)\right| \leq \beta(|\xi|, t) \quad \text{and} \quad \left|y_\lambda(t, \xi, d)\right| \leq \sigma(|h(\xi)|), \qquad \forall\, t \geq 0.$$

Observe that the first inequality together with the UBIBS property implies that system (31) is forward complete, since on any finite interval, one has

$$|x_\lambda(t, \xi, d)| \leq \max\{\sigma(|\xi|), \sigma(\|d\|\,\|y_\lambda\|)\} \leq \max\{\sigma(|\xi|), \sigma(\beta(|\xi|, 0))\}$$

(where we write $x_\lambda(\cdot, \xi, d)$ for the state trajectory associated to the system). Hence, system (31) is UOS. Applying Theorem 2, one sees that there exists some smooth function $V$ such that

- there exist $\alpha_1, \alpha_2 \in \mathcal{K}_\infty$ such that
$$\alpha_1(|h(\xi)|) \leq V(\xi) \leq \alpha_2(|h(\xi)|), \qquad \forall \xi, \tag{32}$$

- there exist some $\alpha_3 \in \mathcal{KL}$ such that
$$DV(\xi)f(\xi, \nu\lambda(|h(\xi)|)) \leq -\alpha_3(V(\xi), |\xi|) \tag{33}$$

for all $\xi \in \mathbb{R}^n$ and all $|\nu| \leq 1$.

It then follows that
$$DV(\xi)f(\xi, \mu) \leq -\alpha_3(V(\xi), |\xi|)$$
whenever $|\mu| \leq \lambda(|h(\xi)|)$, or equivalently, whenever $|h(\xi)| \geq \lambda^{-1}(|\mu|)$. Let $\chi = \alpha_2^{-1} \circ \lambda^{-1}$. Then one has:
$$V(\xi) \geq \chi(|\mu|) \Rightarrow DV(\xi)f(\xi, \mu) \leq -\alpha_3(V(\xi), |\xi|)$$
for all $\xi$ and all $\mu$. Hence, $V$ is an OLIOS-Lyapunov function for the system. ∎

*Sufficiency.* Let $V$ be an OLIOS-Lyapunov function for system (1). Let $\alpha_1, \alpha_2 \in \mathcal{K}_\infty$ such that (10) holds. By (9), and arguing as in Remark 2.1, one also knows that there exists some $\kappa_1$ and $\kappa_2 \in \mathcal{K}_\infty$ such that
$$V(\xi) \geq \chi(|\mu|) \Rightarrow DV(\xi)f(\xi, \mu) \leq -\frac{\kappa_1(V(\xi))}{1 + \kappa_2(|\xi|)} \tag{34}$$
for all $\xi$ and $\mu$.

Let $\beta \in \mathcal{KL}$ be as in Lemma A.4 for the function $\kappa_1$. Pick any initial state $\xi$ and any $u$. Let $x(t)$ and $y(t)$ denote the ensuing trajectory and output function respectively. If for some $t_1 \leq 0$, $V(x(t_1)) \leq \chi(\|u\|)$, then $V(x(t)) \leq \chi(\|u\|)$ for all $t \geq t_1$. (Proof: pick any $\varepsilon > 0$. If $t_1 := \inf\{t > t_1 \,|\, V(x(t)) > \chi(\|u\|) + \varepsilon\}$ is finite, then $V(x(t)) > \chi(\|u\|)$ for all $t$ in some left neighborhood of $t_2$, so $DV(x(t))/dt < 0$ and $V(x(t)) > V(x(t_2))$ for such $t$, contradicting its minimality. As $\varepsilon$ was arbitrary, the claim follows.) Now let
$$\widetilde{t} = \inf\{t \geq 0 : V(x(t)) \leq \chi(\|u\|)\}$$
with the understanding that $\widetilde{t} = \infty$ if $V(x(t)) > \chi(\|u\|)$ for all $t \geq 0$. Then
$$V(x(t)) \leq \chi(\|u\|), \qquad \forall t \geq \widetilde{t}, \tag{35}$$
and on $[0, \widetilde{t})$, it holds that
$$\frac{d}{dt}V(x(t)) \leq -\frac{\kappa_1(V(x(t)))}{1 + \kappa_2(|x(t)|)}.$$
Since the system is UBIBS, there exists some $\sigma$ such that (2) holds. Hence,
$$\frac{d}{dt}V(x(t)) \leq -\frac{\kappa_1(V(x(t)))}{1 + \max\{\widetilde{\kappa}_2(|\xi|), \widetilde{\kappa}_2(\|u\|)\}}$$
for all $t \in [0, \widetilde{t})$, where $\widetilde{\kappa}_2 = \kappa_2 \circ \sigma$. It then follows Lemma A.4 that
$$V(x(t)) \leq \beta\left(V(\xi), \frac{t}{1 + \max\{\widetilde{\kappa}_2(|\xi|), \widetilde{\kappa}_2(\|u\|)\}}\right)$$

for all $t \in [0, \widetilde{t})$.

Let $v_0(s) = \max_{|\xi| \leq s} V(\xi)$. Then $v_0$ is nonincreasing, $v_0(0) = 0$, and $V(\xi) \leq v_0(|\xi|)$. Note then that

$$\beta\left(V(\xi), \frac{t}{1 + \max\{\widetilde{\kappa}_2(|\xi|), \widetilde{\kappa}_2(\|u\|)\}}\right)$$
$$\leq \max\left\{\beta\left(V(\xi), \frac{t}{1 + \widetilde{\kappa}_2(|\xi|)}\right), \beta\left(v_0(\|u\|), \frac{t}{1 + \widetilde{\kappa}_2(\|u\|)}\right)\right\}$$
$$\leq \max\left\{\beta\left(V(\xi), \frac{t}{1 + \widetilde{\kappa}_2(|\xi|)}\right), \beta\left(v_0(\|u\|), 0\right)\right\}$$

(consider two cases: $|\xi| \geq \|u\|$ and $|\xi| \leq \|u\|$). This shows that

$$V(x(t)) \leq \max\left\{\beta\left(V(\xi), \frac{t}{1 + \widetilde{\kappa}_2(|\xi|)}\right), \widetilde{\beta}_0(\|u\|)\right\}$$

for all $t \in [0, \widetilde{t})$, where $\widetilde{\beta}_0(s) = \beta(v_0(s), 0)$. Combining this with (35), one sees that

$$V(x(t)) \leq \max\left\{\beta\left(V(\xi), \frac{t}{1 + \widetilde{\kappa}_2(|\xi|)}\right), \widetilde{\gamma}(\|u\|)\right\} \tag{36}$$

for all $t \geq 0$, where $\widetilde{\gamma}(s) = \widetilde{\beta}(s) + \chi(s)$. Using the fact that $|h(\xi)| \leq \alpha_1^{-1}(V(\xi))$, we conclude that

$$|y(t)| \leq \max\left\{\widetilde{\beta}\left(|h(\xi)|, \frac{t}{1 + \widetilde{\kappa}_2(|\xi|)}\right), \gamma(\|u\|)\right\} \tag{37}$$

for all $t \geq 0$, where $\widetilde{\beta}(s, r) = \alpha_1^{-1}(\beta(\alpha_2(s), r))$, and $\gamma(s) = \alpha_1^{-1}(\widetilde{\gamma}(s))$. ∎

Observe that the proof used in the sufficiency part of Theorem 1 in fact shows that if a system admits an OLIOS-Lyapunov function, then it satisfies an estimate of type (37). Combining this with the necessity part, one concludes as follows:

**Lemma 4.2** A UBIBS system is OLIOS if and only if there exist $\beta \in \mathcal{KL}, \rho \in \mathcal{K}$, and $\gamma \in \mathcal{K}$ such that

$$|y(t, \xi, u)| \leq \beta\left(|h(\xi)|, \frac{t}{1 + \rho(|\xi|)}\right) + \gamma(\|u\|), \qquad \forall t \geq 0,$$

holds for all trajectories of the system. □

### 4.3 Proof of Theorem 1, Part 2.

*Necessity.* Consider an IOS system (1). By Theorem 1 in [27], there exist some locally Lipschitz map $h_0$ and $\chi \in \mathcal{K}_\infty$ with the property that $h_0(\xi) \geq \chi(|h(\xi)|)$ such that the system

$$\dot{x} = f(x, u), \quad y = h_0(x) \tag{38}$$

is OLIOS. By Part 1 of this theorem, system (38) admits an OLIOS-Lyapunov function $V$. This means that there exist $\alpha_1, \alpha_2, \rho \in \mathcal{K}_\infty$ and $\alpha_3 \in \mathcal{KL}$ such that

$$\alpha_1(|h_0(\xi)|) \leq V(\xi) \leq \alpha_2(|h_0(\xi)|), \qquad \forall \xi \in \mathbb{R}^n,$$

and

$$V(\xi) \geq \rho(|\mu|) \implies DV(\xi)f(\xi, \mu) \leq -\alpha_3(V(\xi), |\xi|).$$

To show that $V$ is an IOS-Lyapunov function, it only remains to show that $V(\xi) \geq \widetilde{\alpha}_1(|h(\xi)|)$ for some $\widetilde{\alpha}_1 \in \mathcal{K}_\infty$. But this follows immediately from the fact that $|h(\xi)| \leq \chi^{-1}(h_0(\xi))$. So one can let $\widetilde{\alpha}_1 := \alpha_1 \circ \chi$. Hence, $V$ is indeed an IOS-Lyapunov function for system (1). ∎

*Sufficiency.* Let $V$ be an IOS-Lyapunov function for system (1). From the proof of Part 1 of Theorem 1 (sufficiency), one can see that if $V$ satisfies (34) for some $\chi, \kappa_1, \kappa_2 \in \mathcal{KL}$, then there exist $\widetilde{\beta} \in \mathcal{KL}, \widetilde{\kappa}_2, \widetilde{\gamma} \in \mathcal{K}_\infty$ such that (36) holds. This means that the system

$$\dot{x} = f(x, u), \quad y = V(x)$$

is OLIOS. Since $V(x) \geq \alpha_1(|h(\xi)|)$ for some $\alpha_1 \in \mathcal{K}_\infty$, it follows that system (1) is IOS. ∎

### 4.4 Proof of Theorem 1, Part 3.

*Necessity.* Since the system (1) is ROS, there is a smooth $\mathcal{K}_\infty$-function $\lambda$ such that (6) holds for the corresponding system (7). As remarked in the necessity proof of Theorem 1, part 1, one sees that system (7) is also forward complete, and hence, UOS. By Theorem 2, system (7) admits a smooth Lyapunov function $V$ satisfying (25) and

$$DV(\xi)f(\xi, \mu\lambda(|y|)) - \alpha_3(V(\xi), |\xi|), \qquad \forall \xi \in \mathbb{R}^n, \ \forall \ |\mu| \leq 1,$$

for some $\alpha_3 \in \mathcal{KL}$. This is equivalent to

$$|y| \geq \lambda^{-1}(|\nu|) \ \Rightarrow \ DV(\xi)f(\xi, \nu) - \alpha_3(V(\xi), |\xi|), \qquad \forall \xi \in \mathbb{R}^n, \ \forall \ |\nu| \in \mathbb{R}^m.$$

Hence, one concludes that $V$ is an ROS-Lyapunov function for system (1).

*Sufficiency.* Let $V$ be an ROS-Lyapunov function. As in Remark 2.1, there exist $\chi, \kappa_1, \kappa_2 \in \mathcal{K}_\infty$ such that

$$DV(\xi)f(\xi, \mu) \leq -\frac{\kappa_1(V(\xi))}{1 + \kappa_2(|\xi|)}$$

whenever $|h(\xi)| \geq \chi(|\mu|)$. Let $\lambda = \chi^{-1}$. Without loss of generality, one may assume that $\lambda$ is smooth. (Otherwise, one can always replace $\lambda$ by a smooth $\mathcal{K}_\infty$-function that is majorized by $\lambda$.) It then follows that

$$DV(\xi)f(\xi, \nu\lambda(|h(\xi)|)) \leq -\frac{\kappa_1(V(\xi))}{1 + \kappa_2(|\xi|)}$$

for all $\xi \in \mathbb{R}^n$ and all $|\nu| \leq 1$. This implies that for any trajectory $x_\lambda(t) = x_\lambda(t, \xi, d)$ of the system

$$\dot{x} = f(x, d\lambda(|y|)), \quad y = h(x),$$

where $d \in \mathcal{M}_\mathcal{B}$, it holds that

$$\frac{d}{dt}V(x_\lambda(t)) \leq -\frac{\kappa_1(V(x_\lambda(t)))}{1 + \kappa_2(|x_\lambda(t)|)} \tag{39}$$

for all $t \geq 0$. It follows immediately that $V(x_\lambda(t)) \leq V(\xi)$ for all $t \geq 0$. Since $V(\xi) \geq \alpha_1(|h(\xi)|)$ for some $\alpha_1 \in \mathcal{K}_\infty$, it follows that, for some $\sigma \in \mathcal{K}_\infty$,

$$|y_\lambda(t)| \leq \sigma(|\xi|), \qquad \forall \, t \geq 0. \tag{40}$$

Since the system is UBIBS, there exists some $\sigma_0 \in \mathcal{K}$ such that

$$|x_\lambda(t, \xi, d)| \leq \max\{\sigma_0(|\xi|), \sigma_0(\|u_d\|)\}, \qquad \forall \, t \geq 0,$$

where $u_d(t) = d(t)\lambda(|y(t)|)$. Combining this with (40), it follows that

$$\left|x_\lambda(t,\xi,d)\right| \leq \widetilde{\sigma}(|\xi|), \qquad \forall\, t \geq 0,$$

where $\widetilde{\sigma}(s) = \max\{\sigma_0(s), \sigma_0(\lambda(\sigma(s)))\}$. Substituting this back into (39), one has

$$\frac{d}{dt} V(x_\lambda(t)) \leq -\frac{\kappa_1(V(x_\lambda(t)))}{1 + \kappa_3(|\xi|)}, \qquad \forall\, t \geq 0,$$

where $\kappa_3(s) = \kappa_2(\widetilde{\sigma}(s))$. Again, by Lemma A.4, one knows that there exists some $\beta \in \mathcal{KL}$ (which depends only upon $\kappa_1$) such that

$$V(x_\lambda(t)) \leq \beta\left(V(\xi), \frac{t}{1 + \kappa_3(|\xi|)}\right), \qquad \forall\, t \geq 0.$$

Together with the fact that $|h(\xi)| \leq \alpha_1^{-1}(V(\xi))$, this yields

$$\left|y_\lambda(t,\xi,d)\right| \leq \widetilde{\beta}(|\xi|, t), \qquad \forall\, t \geq 0,$$

where $\widetilde{\beta}(s,r) = \alpha_1^{-1}[\beta(\alpha_2(s), t/(1 + \kappa_3(s)))]$ is in $\mathcal{KL}$, and $\alpha_2$ is any $\mathcal{K}_\infty$-function such that $V(\xi) \leq \alpha_2(|\xi|)$ for all $\xi$. This shows that the system is ROS. ∎

### 4.5 Proof of Theorem 1, Part 4.

*Necessity.* Assume that a UBIBS system (1) admits an estimate (5) for some $\beta \in \mathcal{KL}$ and some $\gamma \in \mathcal{K}$. Then it admits an estimate of type (4) with $\sigma_1(s) = 2\beta(s,0)$ and $\sigma_2(s) = 2\gamma(s)$. By following the proof of Lemma 4.1, one can show that there exists some $\lambda \in \mathcal{K}_\infty$ such that for any $\xi$ and $u$, if $\gamma(|u(t)|) \leq \lambda(|y(t,\xi,u)|)$ for almost all $t \geq 0$, it holds that $\gamma(|u(t)|) \leq \frac{|h(\xi)|}{2}$. Again, without loss of generality, one may assume that $\lambda$ is smooth. One then can show that for the system

$$\dot{x}(t) = f(x(t), d(t)\lambda(|y(t)|)), \quad y(t) = h(x(t)),$$

there exists $\widetilde{\beta} \in \mathcal{KL}$ so that, for all trajectories $x_\lambda(t,\xi,d)$, it holds that

$$\left|y_\lambda(t,\xi,d)\right| \leq \widetilde{\beta}(|h(\xi)|, t)$$

for all $t \geq 0$. Applying the last part of Theorem 2 one sees that there exists $V$ satisfying (27) for some $\alpha_1, \alpha_2 \in \mathcal{K}_\infty$ and

$$DV(\xi)f(\xi, \nu\lambda(|y(\xi)|)) \leq -\alpha_3(V(\xi))$$

for all $\xi$ and all $|\nu| \leq 1$. This is equivalent to the existence of $\chi \in \mathcal{K}_\infty$ such that

$$V(\xi) \geq \chi(|\mu|) \Rightarrow DV(\xi)f(\xi, u) \leq -\alpha_3(V(\xi)). \tag{41}$$

*Sufficiency.* It is routine to show that if there is a smooth function $V$ satisfying (27) and (41), then the system admits an estimate of type (5). ∎

## 4.6 Proof of Theorem 2

Consider the system
$$\dot{x}(t) = f(x(t), u(t)), \quad y = h(x(t)), \tag{42}$$
where the input $u$ takes values in a compact subset $\Omega$ of $\mathbb{R}^m$. Assume that the system is UBIBS and there exists some $\beta \in \mathcal{KL}$ such that (22) holds for all trajectories of (42). Let $\omega : \mathbb{R}^n \to \mathbb{R}_{\geq 0}$ be defined by
$$\omega(\xi) := \sup\{|y(t, \xi, u)| : t \geq 0, u \in \mathcal{M}_\Omega\}. \tag{43}$$

It then holds that
$$|h(\xi)| \leq \omega(\xi) \leq \beta_0(|\xi|), \quad \forall \xi \in \mathbb{R}^n, \tag{44}$$
where $\beta_0(s) = \beta(s, 0)$. Moverover, if there exists some $\sigma \in \mathcal{K}$ such that (23) holds for all trajectories, then the above can be strengthened to
$$|h(\xi)| \leq \omega(\xi) \leq \sigma(|h(\xi)|), \quad \forall \xi \in \mathbb{R}^n. \tag{45}$$

Observe that, for any $\xi \in \mathbb{R}^n$, $u \in \mathcal{M}_\Omega$, and $t_1 \geq 0$,
$$\omega(x(t_1, \xi, u)) \leq \sup_{t \geq 0, v \in \mathcal{M}_\Omega} |y(t_1 + t, \xi, v)| \leq \beta(|\xi|, t_1). \tag{46}$$

Also $\omega$ decreases along trajectories, i.e.,
$$\omega(x(t, \xi, u)) \leq \omega(\xi), \quad \forall t \geq 0, \, \xi \in \mathbb{R}^n, \, u \in \mathcal{M}_\Omega. \tag{47}$$

Define
$$\mathcal{D} := \{\xi : y(t, \xi, u) = 0, \, \forall t \geq 0, \, \forall u \in \mathcal{M}_\Omega\}.$$

Then $\omega(\xi) = 0$ if and only if $\xi \in \mathcal{D}$. For $\xi \notin \mathcal{D}$, it holds that
$$\omega(\xi) = \sup_{0 \leq t \leq t_\xi, \, u \in \mathcal{M}_\Omega} |y(t, \xi, u)|, \tag{48}$$
where $t_\xi = T_{|\xi|}(\omega(\xi)/2)$, and $T_r(s)$ is defined as in Lemma A.1 associated with the function $\beta$.

**Lemma 4.3** The function $\omega(\xi)$ is locally Lipschitz on $\mathbb{R}^n \setminus \mathcal{D}$ and continuous everywhere.

*Proof.* First notice that
$$\varliminf_{\xi \to \xi_0} \omega(\xi) \geq \omega(\xi_0), \quad \forall \xi_0 \in \mathbb{R}^n, \tag{49}$$
that is, $\omega(\xi)$ is lower semi-continuous on $\mathbb{R}^n$. Indeed, pick $\xi_0$ and let $c := \omega(\xi_0)$. Take any $\varepsilon > 0$. Then there are some $u_0$ and $t_0$ so that $|y(t_0, \xi_0, u_0)| \geq c - \varepsilon/2$. By continuity of $y(t_0, \cdot, u_0)$, there is some neighborhood $\widetilde{U}_0$ of $\xi_0$ so that $|y(t_0, \xi, u_0)| \geq c - \varepsilon$ for all $\xi \in \widetilde{U}_0$. Thus $\omega(\xi) \geq c - \varepsilon$ for all $\xi \in \widetilde{U}_0$, and this establishes (49).

Fix any $\xi_0 \in \mathbb{R}^n \setminus \mathcal{D}$, and let $c_0 = \omega(\xi_0)/2$. Then there exists a bounded neighborhood $U_0$ of $\xi_0$ such that
$$\omega(\xi) \geq c_0, \quad \forall \xi \in U_0.$$

Let $s_0$ be such that $|\xi| \leq s_0$ for all $\xi \in U_0$. Then
$$\omega(\xi) = \sup\{|y(t, \xi, u)| : t \in [0, t_1], u \in \mathcal{M}_\Omega\}, \quad \forall \xi \in U_0,$$

where $t_1 = T_{s_0}(c_0/2)$. By [13, Proposition 5.5], one knows that $x(t, \xi, u)$ is locally Lipschitz in $\xi$ uniformly on $u \in \mathcal{M}_\Omega$ and on $t \in [0, t_1]$, and therefore, so is $y(t, \xi, u)$. Let $C$ be a constant such that

$$|y(t, \xi, u) - y(t, \eta, u)| \leq C |\xi - \eta|, \qquad \forall \xi, \eta \in U_0, \ \forall 0 \leq t \leq t_1, \ \forall u \in \mathcal{M}_\Omega.$$

For any $\varepsilon > 0$ and any $\xi \in U_0$, there exist some $t_{\xi,\varepsilon} \in [0, t_1]$ and some $u_{\xi,\varepsilon}$ such that

$$\omega(\xi) \leq |y(t_{\xi,\varepsilon}, \xi, u_{\xi,\varepsilon})| + \varepsilon.$$

It then follows that, for any $\xi, \eta \in U_0$, for any $\varepsilon > 0$,

$$\omega(\xi) - \omega(\eta) \leq |y(t_{\xi,\varepsilon}, \xi, u_{\xi,\varepsilon})| + \varepsilon - |y(t_{\xi,\varepsilon}, \eta, u_{\xi,\varepsilon})| \leq C |\xi - \eta| + \varepsilon.$$

Consequently,

$$\omega(\xi) - \omega(\eta) \leq C |\xi - \eta|, \qquad \forall \xi, \eta \in U_0.$$

By symmetry,

$$\omega(\eta) - \omega(\xi) \leq C |\xi - \eta|, \qquad \forall \xi, \eta \in U_0.$$

This proves that $\omega$ is locally Lipschitz on $\mathbb{R}^n \setminus \mathcal{D}$.

We now show that $\omega$ is continuous on $\mathcal{D}$. Fix $\xi_0 \in \mathcal{D}$. One would like to show that

$$\lim_{\xi \to \xi_0} \omega(\xi) = 0. \tag{50}$$

Assume that this does not hold. Then there exists some $\varepsilon_0 > 0$ and a sequence $\{\xi_k\}$ with $\xi_k \to \xi_0$ such that $\omega(\xi_k) > \varepsilon_0$ for all $k$. Without loss of generality, one may assume that

$$|\xi_k| \leq s_1, \qquad \forall k,$$

for some $s_1 \geq 0$. It then follows that

$$\omega(\xi_k) = \sup \{|y(t, \xi_k, u)| : \ t \in [0, t_2], u \in \mathcal{M}_\Omega\},$$

where $t_2 = T_{s_1}(\varepsilon_0/2)$. Hence, for each $k$, there exists some $u_k \in \mathcal{M}_\Omega$ and some $\tau_k \in [0, t_2]$ such that

$$|y(\tau_k, \xi_k, u_k)| \geq \omega(\xi_k) - \varepsilon_0/2 \geq \varepsilon_0/2.$$

Again, by the locally Lipschitz continuity of the trajectories, one knows that there is some $C_1 > 0$ such that

$$|y(t, \xi_k, u) - y(t, \xi_0, u)| \leq C_1 |\xi_k - \xi_0|, \qquad \forall k \geq 0, \ \forall 0 \leq t \leq t_2, \ \forall u \in \mathcal{M}_\Omega.$$

Hence,

$$|y(\tau_k, \xi_0, u_k)| \geq \varepsilon_0/4$$

for $k$ large enough, contradicting the fact that $y(t, \xi_0, u) \equiv 0$ for all $u \in \mathcal{M}_\Omega$. This shows that (50) holds if $\xi_0 \in \mathcal{D}$. ∎

Next, we pick any smooth and bounded function $k : \mathbb{R}_{\geq 0} \to \mathbb{R}_{>0}$ whose derivative is everywhere positive, and define $W : \mathbb{R}^n \to \mathbb{R}_{\geq 0}$ by

$$W(\xi) := \sup \{\omega(x(t, \xi, u))k(t) : \ t \geq 0, \ u \in \mathcal{M}_\Omega\}. \tag{51}$$

Corresponding to $k$ there are two positive real numbers $c_1 < c_2$ such that $k(t) \in [c_1, c_2]$ for all $t \geq 0$, and so
$$c_1 \omega(\xi) \leq W(\xi) \leq c_2 \omega(\xi), \qquad \forall \xi \in \mathbb{R}^n,$$
which implies that
$$c_1 |h(\xi)| \leq W(\xi) \leq c_2 \beta_0(|\xi|), \qquad \forall \xi \in \mathbb{R}^n. \tag{52}$$

Note, for future reference, that it is always possible to find a bounded, positive, and decreasing continuous function $\tau(\cdot)$ with $\tau(t) \to 0$ as $t \to \infty$, such that
$$k'(t) \geq \tau(t) \quad \text{for all } t \geq 0. \tag{53}$$

By (46), one knows that $\omega(x(t, \xi, u)) \to 0$ as $t \to \infty$. It follows that there is some $\tau_\xi \geq 0$ such that
$$W(\xi) = \sup \{\omega(x(t, \xi, u)) k(t) : u \in \mathcal{M}_\Omega, \ 0 \leq t \leq \tau_\xi\}. \tag{54}$$
Furthermore, one can get the following estimate, where $\{T_r\}$ is a family of functions associated to $\beta$ as in Lemma A.4.

**Lemma 4.4** For any $\xi \notin \mathcal{D}$ with $|\xi| \leq r$,
$$W(\xi) = \sup \{\omega(x(t, \xi, u)) k(t) : u \in \mathcal{M}_\Omega, \ 0 \leq t \leq \tau_\xi\},$$
where $\tau_\xi = T_r(\frac{c_1}{2c_2} \omega(\xi))$.

*Proof.* If the statement is not true, then for any $\varepsilon > 0$, there exists some $t_\varepsilon > T_r(\frac{c_1}{2c_2} \omega(\xi))$ and some $u_\varepsilon \in \mathcal{M}_\Omega$ such that
$$W(\xi) \leq \omega(x(t_\varepsilon, \xi, u_\varepsilon)) k(t_\varepsilon) + \varepsilon.$$
This implies the following:
$$\begin{aligned}
\omega(\xi) &\leq \frac{1}{c_1} W(\xi) \leq \frac{1}{c_1} \omega(x(t_\varepsilon, \xi, u_\varepsilon)) k(t_\varepsilon) + \frac{\varepsilon}{c_1} \\
&\leq \frac{c_2}{c_1} \omega(x(t_\varepsilon, \xi, u_\varepsilon)) + \frac{\varepsilon}{c_1} \leq \frac{c_2}{c_1} \cdot \frac{c_1}{2c_2} \omega(\xi) + \frac{\varepsilon}{c_1} \\
&= \frac{\omega(\xi)}{2} + \frac{\varepsilon}{c_1}.
\end{aligned}$$
Taking the limit as $\varepsilon \to 0$ results in a contradiction.

**Lemma 4.5** The function $W(\cdot)$ is locally Lipschitz on $\mathbb{R}^n \setminus \mathcal{D}$ and continuous everywhere.

*Proof.* Fix $\xi_0 \notin \mathcal{D}$. Let $K_0$ be a compact neighborhood of $\xi_0$ such that $K_0 \cap \mathcal{D} = \emptyset$. Since $\omega$ is continuous, it follows that there is some $r_0 > 0$ such that $\omega(\xi) > r_0$ for all $\xi \in K_0$, and hence, $W(\xi) > r_1 := c_1 r_0$ for all $\xi \in K_0$. Let
$$T_0 = T_{s_0}\left(\frac{r_1}{8c_2}\right),$$
where $s_0 > 0$ is such that $|\xi| \leq s_0$ for all $\xi \in K_0$. Let $C > 0$ be such that
$$|y(t, \xi, u) - y(t, \eta, u)| \leq C |\xi - \eta|, \quad \forall t \in [0, T_0], \ \forall \xi, \eta \in K_0, \ \forall u \in \mathcal{M}_\Omega.$$

Let
$$K_1 = K_0 \cap \left\{\xi : |\xi - \xi_0| \leq \frac{r_1}{16Cc_2}\right\}.$$

Fix any $\varepsilon \in (0, r_1/4)$. Then, for any $\xi \in K_1$, there exist $t_{\xi,\varepsilon} \in [0, T_0]$ and $u_{\xi,\varepsilon} \in \mathcal{M}_\Omega$ such that
$$W(\xi) \leq \omega(x(t_{\xi,\varepsilon}, \xi, u_{\xi,\varepsilon}))k(t_{\xi,\varepsilon}) + \varepsilon.$$

*Claim:* For any $\xi, \eta \in K_1$, $\omega(x(t_{\xi,\varepsilon}, \eta, u_{\xi,\varepsilon})) \geq \frac{r_1}{8c_2}$.

*Proof.* First we note that for any $\xi \in K_1 \subset K_0$,
$$\omega(x(t_{\xi,\varepsilon}, \xi, u_{\xi,\varepsilon})) \geq \frac{W(\xi) - \varepsilon}{c_2} \geq \frac{W(\xi)}{2c_2} \geq r_2,$$

where $r_2 := \frac{r_1}{2c_2}$. Thus, for each $\xi \in K_1$, there exists some $v_\xi \in \mathcal{M}_\Omega$ and some $\tau_\xi > 0$ such that
$$|y(\tau_\xi, x(t_{\xi,\varepsilon}, \xi, u_{\xi,\varepsilon}), v_\xi)| \geq \omega(x(t_{\xi,\varepsilon}, \xi, u_{\xi,\varepsilon})) - r_2/2 \geq r_2/2.$$

Observe that
$$y(\tau_\xi, x(t_{\xi,\varepsilon}, \xi, u_{\xi,\varepsilon}), v_\xi) = y(\tau_\xi + t_{\xi,\varepsilon}, \xi, \overline{v}_{\xi,\varepsilon}),$$

where $\overline{v}_{\xi,\varepsilon}$ is the concatenation of $u_{\xi,\varepsilon}$ and $v_\xi$, i.e.,
$$\overline{v}_{\xi,\varepsilon}(t) = \begin{cases} u_{\xi,\varepsilon}(t), & \text{if } 0 \leq t < t_{\xi,\varepsilon}, \\ v_\xi(t - t_{\xi,\varepsilon}), & \text{if } t \geq t_{\xi,\varepsilon}. \end{cases}$$

Noticing that $|y(t, \xi, u)| \leq r_2/2$ for all $t \geq T_{s_0}(r_2/4)$, one concludes that $\tau_\xi + t_{\xi,\varepsilon} < T_{s_0}(r_2/4) = T_0$. Note also that for any $\eta \in K_1$,
$$\begin{aligned}
|y(\tau_\xi, x(t_{\xi,\varepsilon}, \eta, u_{\xi,\varepsilon}), v_\xi)| &= |y(\tau_\xi + t_{\xi,\varepsilon}, \eta, \overline{v}_{\xi,\varepsilon})| \\
&\geq |y(\tau_\xi + t_{\xi,\varepsilon}, \xi, \overline{v}_{\xi,\varepsilon})| - |y(\tau_\xi + t_{\xi,\varepsilon}, \eta, \overline{v}_{\xi,\varepsilon}) - y(\tau_\xi + t_{\xi,\varepsilon}, \xi, \overline{v}_{\xi,\varepsilon})| \\
&\geq \frac{r_2}{2} - C|\xi - \eta| \\
&\geq \frac{r_2}{2} - 2C\frac{r_1}{16Cc_2} = \frac{r_1}{4c_2} - \frac{r_1}{8c_2} = \frac{r_1}{8c_2}.
\end{aligned}$$

This implies that $\omega(x(t_{\xi,\varepsilon}, \eta, u_{\xi,\varepsilon})) \geq \frac{r_1}{8c_2}$ for all $\xi, \eta \in K_1$, as claimed.

According to [13, Proposition 5.1], there is some compact set $K_2$ such that $x(t, \xi, u) \in K_2$ for all $0 \leq t \leq T_0$, all $\xi \in K_1$ and all $u \in \mathcal{M}_\Omega$. Let
$$K_3 = K_2 \cap \{\xi : \omega(\xi) \geq r_1/8c_2\}.$$

Applying Lemma 4.3, one knows that there is some $C_1 > 0$ such that
$$|\omega(\zeta_1) - \omega(\zeta_2)| \leq C_1 |\zeta_1 - \zeta_2|, \qquad \forall \zeta_1, \zeta_2 \in K_3.$$

Since for all $\xi, \eta \in K_1$ and all $0 < \varepsilon < r_1/4$, $x(t_{\xi,\varepsilon}, \eta, u_{\xi,\varepsilon}) \in K_3$, we have:
$$|\omega(x(t_{\xi,\varepsilon}, \xi, u_{\xi,\varepsilon})) - \omega(x(t_{\xi,\varepsilon}, \eta, u_{\xi,\varepsilon}))| \leq C_1 |x(t_{\xi,\varepsilon}, \xi, u_{\xi,\varepsilon}) - x(t_{\xi,\varepsilon}, \eta, u_{\xi,\varepsilon})|,$$

for all $\xi, \eta \in K_1$, all $\varepsilon \in (0, r_1/4)$. Hence,
$$\begin{aligned}
W(\xi) - W(\eta) &\leq \omega(x(t_{\xi,\varepsilon}, \xi, u_{\xi,\varepsilon}))k(t_{\xi,\varepsilon}) - \omega(x(t_{\xi,\varepsilon}, \eta, u_{\xi,\varepsilon}))k(t_{\xi,\varepsilon}) + \varepsilon \\
&\leq c_2 |\omega(x(t_{\xi,\varepsilon}, \xi, u_{\xi,\varepsilon})) - \omega(x(t_{\xi,\varepsilon}, \eta, u_{\xi,\varepsilon}))| + \varepsilon \\
&\leq c_2 C_1 |x(t_{\xi,\varepsilon}, \xi, u_{\xi,\varepsilon}) - x(t_{\xi,\varepsilon}, \eta, u_{\xi,\varepsilon})| + \varepsilon \\
&\leq c_2 C_1 C_2 |\xi - \eta| + \varepsilon,
\end{aligned}$$

where $C_2 > 0$ is such a constant that $|x(t, \xi, u) - x(t, \eta, u)| \leq C_2 |\xi - \eta|$ for all $\xi, \eta \in K_3$, all $t \in [0, T_0]$, all $u \in \mathcal{M}_\Omega$. Note that the above holds for any $\varepsilon \in (0, r_1/4)$, and thus,

$$W(\xi) - W(\eta) \leq C_3 |\xi - \eta|$$

for all $\xi, \eta \in K_1$, where $C_3 = c_2 C_1 C_2$. By symmetry, one proves that

$$W(\eta) - W(\xi) \leq C_3 |\xi - \eta|$$

for all $\xi, \eta \in K_1$.

To prove the continuity of $W$ on $\mathcal{D}$, it is enough to notice that for any $\xi \in \mathcal{D}$, $W(\xi) = 0$ and

$$|W(\xi) - W(\eta)| \leq c_2 \omega(\eta) \to 0, \qquad \text{as } \eta \to \xi.$$

The proof of Lemma 4.5 is thus concluded. ∎

Below we show that $W$ is decreasing along trajectories. Pick any $\xi \notin \mathcal{D}$. Let $\theta_0 > 0$ be such that

$$\omega(x(t, \xi, \mathbf{v})) \geq \omega(\xi)/2, \quad \forall t \in [0, \theta_0], \ \forall v \in \Omega,$$

where $\mathbf{v}$ denotes the constant function $\mathbf{v}(t) \equiv v$. (Observe that such a $\theta_0$ exists because $\omega$ is continuous.) Pick any $\theta \in [0, \theta_0]$, and let $\eta_\mathbf{v} = x(\theta, \xi, \mathbf{v})$. For any $\varepsilon > 0$, there exists some $t_{\mathbf{v}, \varepsilon}$ and $u_{\mathbf{v}, \varepsilon} \in \mathcal{M}_\Omega$ such that

$$\begin{aligned}
W(\eta_\mathbf{v}) &\leq \omega(x(t_{\mathbf{v},\varepsilon}, \eta_\mathbf{v}, u_{\mathbf{v},\varepsilon})) k(t_{\mathbf{v},\varepsilon}) + \varepsilon \\
&= \omega(x(t_{\mathbf{v},\varepsilon} + \theta, \xi, \overline{u}_{\mathbf{v},\varepsilon})) k(t_{\mathbf{v},\varepsilon} + \theta) \left(1 - \frac{k(t_{\mathbf{v},\varepsilon} + \theta) - k(t_{\mathbf{v},\varepsilon})}{k(t_{\mathbf{v},\varepsilon} + \theta)}\right) + \varepsilon \\
&\leq W(\xi) \left(1 - \frac{k(t_{\mathbf{v},\varepsilon} + \theta) - k(t_{\mathbf{v},\varepsilon})}{c_2}\right) + \varepsilon,
\end{aligned} \tag{55}$$

where $\overline{u}_{\mathbf{v},\varepsilon}$ denotes the concatenation of $\mathbf{v}$ and $u_{\mathbf{v},\varepsilon}$. Still for the fixed $\xi$ and $\theta$, and for any $r > |\xi|$, define

$$T^r_{\xi,\theta} := \max_{\widetilde{v} \in \Omega} T_r\left(\frac{c_1}{2c_2} \omega(x(\theta, \xi, \widetilde{\mathbf{v}}))\right). \tag{56}$$

Notice that $x(\theta, \xi, \widetilde{\mathbf{v}})$ is jointly continuous as a function of $(\theta, \xi, \widetilde{v})$. Since $\omega$ and $T_r$ are both continuous, this maximum is well defined and, moreover, $T^r_{\xi,\theta}$ is continuous as a function of $\theta$, so in particular

$$\lim_{\theta \to 0^+} T^r_{\xi,\theta} = T_r\left(\frac{c_1}{c_2} \omega(\xi)\right). \tag{57}$$

*Claim:* $t_{\mathbf{v},\varepsilon} + \theta \leq T^r_{\xi,\theta}$ for all $v \in \Omega$ and for all $\varepsilon \in (0, \frac{c_1}{4}\omega(\xi))$.

*Proof.* Assume that this is not true. Then there is some $v \in \Omega$ and some $\varepsilon \in \left(0, \frac{c_1}{4}\omega(\xi)\right)$ such that $t_{\mathbf{v},\varepsilon} + \theta > T^r_{\xi,\theta}$, and in particular,

$$t_{\mathbf{v},\varepsilon} + \theta \geq T_r\left(\frac{c_1}{2c_2} \omega(x(\theta, \xi, \mathbf{v}))\right),$$

from which it follows that

$$\omega(x(t_{\mathbf{v},\varepsilon}, \eta_\mathbf{v}, u_{\mathbf{v},\varepsilon})) = \omega(x(t_{\mathbf{v},\varepsilon} + \theta, \xi, \overline{u}_{\mathbf{v},\varepsilon})) \leq \frac{c_1}{2c_2} \omega(x(\theta, \xi, \mathbf{v})) = \frac{c_1}{2c_2} \omega(\eta_\mathbf{v}),$$

for some input function $\overline{u}_{\mathbf{v},\varepsilon}$ (which we can take to be the concatenation of $\mathbf{v}$ and $u_{\mathbf{v},\varepsilon}$; note that the inequality follows from (46) and the definition of the functions $T_r$).

By the definition of $W$, one has

$$\begin{aligned}
\omega(\eta_{\mathbf{v}}) &\leq \frac{1}{c_1} W(\eta_{\mathbf{v}}) \leq \frac{1}{c_1} \omega(t_{\mathbf{v},\varepsilon}, \eta_{\mathbf{v}}, u_{\mathbf{v},\varepsilon}) k(t_{\mathbf{v},\varepsilon}) + \frac{\varepsilon}{c_1} \\
&\leq \frac{c_2}{c_1} \omega(t_{\mathbf{v},\varepsilon} + \theta, \xi, \overline{u}_{\mathbf{v},\varepsilon}) + \frac{\varepsilon}{c_1} \\
&\leq \frac{1}{2} \omega(\eta_{\mathbf{v}}) + \frac{\varepsilon}{c_1},
\end{aligned}$$

which is impossible, since $\varepsilon < \frac{c_1}{4}\omega(\xi) \leq \frac{c_1}{2}\omega(\eta_{\mathbf{v}})$. This proves the claim.

From (55), we have, for any $v \in \mathcal{D}$ and for any $\varepsilon$ small enough,

$$W(x(\theta, \xi, \mathbf{v})) - W(\xi) \leq -\frac{W(\xi)}{c_2} \tau(t_{\mathbf{v},\varepsilon} + c\theta)\theta + \varepsilon$$

for some $c \in (0, 1)$, where we used the mean value theorem in order to estimate the change in $k$, and where $\tau$ is a function as in (53). Using the monotonicity of $\tau(\cdot)$, and the above claim, one concludes

$$W(x(\theta, \xi, \mathbf{v})) - W(\xi) \leq -\frac{W(\xi)}{c_2} \tau\left(T^r_{\xi,\theta}\right) \theta + \varepsilon$$

for all $\varepsilon$ small enough. Letting $\varepsilon \to 0$, one obtains:

$$W(x(\theta, \xi, \mathbf{v})) - W(\xi) \leq -\frac{W(\xi)}{c_2} \tau\left(T^r_{\xi,\theta}\right) \theta, \quad \forall v \in \Omega.$$

Thus one concludes that for any $v \in \Omega$ and any $\theta > 0$,

$$\frac{W(x(\theta, \xi, \mathbf{v})) - W(\xi)}{\theta} \leq -\frac{W(\xi)}{c_2} \tau(T^r_{\xi,\theta}).$$

Since $W$ is locally Lipschitz on $\mathbb{R}^n \setminus \mathcal{D}$, it is differentiable almost everywhere on $\mathbb{R}^n \setminus \mathcal{D}$, and hence, for any $v \in \Omega$, any $r > |\xi|$, and any $\xi$ at which $W$ is differentiable,

$$\begin{aligned}
DW(\xi)f(\xi, v) &= \lim_{\theta \to 0^+} \frac{W(x(\theta, \xi, \mathbf{v})) - W(\xi)}{\theta} \leq -\lim_{\theta \to 0^+} \frac{W(\xi)}{c_2} \tau(T^r_{\xi,\theta}) \\
&= -\frac{W(\xi)}{c_2} \tau\left(\lim_{\theta \to 0^+} T^r_{\xi,\theta}\right) = -\frac{W(\xi)}{c_2} \tau\left(T_r\left(\frac{c_1}{c_2}\omega(\xi)\right)\right) \\
&\leq -\frac{W(\xi)}{c_2} \tau\left(T_r\left(\frac{c_1}{c_2^2}W(\xi)\right)\right) = -\widetilde{\alpha}_3(W(\xi), r),
\end{aligned} \qquad (58)$$

where $\widetilde{\alpha}_3(s, r) = \frac{s}{c_2}\tau(T_r(c_3 s))$ with $c_3 = c_1/c_2^2$. Since (58) holds for all $r > |\xi|$, it follows that

$$DW(\xi)f(\xi, v) \leq -\widetilde{\alpha}_3(W(\xi), 2|\xi|), \qquad (59)$$

for all $v \in \Omega$ and for almost all $\xi \in \mathbb{R}^n \setminus \mathcal{D}$.

Since $T_r(s)$ is defined for all $r \geq 0$ and $s > 0$, one sees that $\widetilde{\alpha}_3$ is defined on $\mathbb{R}_{>0} \times \mathbb{R}_{\geq 0}$. Extend $\widetilde{\alpha}_3$ to $\mathbb{R}_{\geq 0} \times \mathbb{R}_{\geq 0}$ by letting $\widetilde{\alpha}_3(0, r) := 0$ for all $r \geq 0$. By the continuity property of $\tau$ and $T_r(\cdot)$, one sees that $\widetilde{\alpha}_3(\cdot, r)$ is continuous for each $r$. (The continuity at $s = 0$ follows from $\widetilde{\alpha}_3(s, r) = s\tau(T_r(c_3 s))/c_2 \leq s\tau(0)/c_2$ for all $s > 0$.) Furthermore, since $\tau(T_r(c_3 s))$ is nondecreasing in $s$, it follows that $\widetilde{\alpha}_3(s,r)$ is of class $\mathcal{K}$ in $s$. Let $\check{\alpha}_3(s,r) = \widetilde{\alpha}_3(s,2r)/(1+r)$. This function tends to zero as $r \to \infty$, because $\widetilde{\alpha}_3(s,r)$ is nonincreasing in $r$; thus $\check{\alpha}_3(s,r)$ is of class $\mathcal{KL}$. Moreover,

$$DW(\xi)f(\xi,v) \leq -\check{\alpha}_3(V(\xi),|\xi|), \qquad \forall \xi \in \mathbb{R}^n \setminus \mathcal{D}, \ \forall v \in \Omega.$$

By Corollary A.3, there exists a continuous $\mathcal{KL}$-function $\widehat{\alpha}_3$ such that

$$DW(\xi)f(\xi,v) \leq -\widehat{\alpha}_3(V(\xi),|\xi|), \qquad \forall \xi \in \mathbb{R}^n \setminus \mathcal{D}, \ \forall v \in \Omega. \tag{60}$$

To complete the proof, we follow the strategy used in [13] to find a smooth approximation of $W$. First of all, by Theorem B.1 in [13], applied on $\mathbb{R}^n \setminus \mathcal{D}$, there is a continuous function $W_1$ that is smooth on $\mathbb{R}^n \setminus \mathcal{D}$ such that

$$|W_1(\xi) - W(\xi)| \leq \frac{W(\xi)}{2}, \qquad \forall \xi \in \mathbb{R}^n \setminus \mathcal{D}, \tag{61}$$

and

$$DW_1(\xi)f(\xi,v) \leq -\widehat{\alpha}_3(W(\xi),|\xi|)/2, \qquad \forall \xi \in \mathbb{R}^n \setminus \mathcal{D}, \ \forall v \in \Omega. \tag{62}$$

We extend $W_1$ to all of $\mathbb{R}^n$ by letting $W_1 \equiv 0$ on $\mathcal{D}$; thus, the approximation (61) holds on all of $\mathbb{R}^n$. (Note that $W$ and $\widehat{\alpha}_3(V(\xi),|\xi|)$ are both continuous, so the result in [13] can indeed be applied.)

Next, we appeal to Lemma 4.3 in [13]. This shows that there exists some $\rho \in \mathcal{K}_\infty$ with $\rho'(s) > 0$ for all $s > 0$ such that $\rho \circ W_1$ is smooth everywhere. Let $V = \rho \circ W_1$. It follows from (52) and (61) that

$$\alpha_1(|h(\xi)|) \leq V(\xi) \leq \alpha_2(|\xi|), \qquad \forall \xi \in \mathbb{R}^n,$$

where $\alpha_1(s) = \rho(c_1 s/2)$, $\alpha_2(s) = \rho(2c_2\beta_0(s))$, and it follows from (61) and (62) that

$$DV(\xi)f(\xi,\mu) \leq -\rho'(W_1(\xi))\widehat{\alpha}_3(W(\xi),|\xi|)/2 \leq -\alpha_3(V(\xi),|\xi|) \tag{63}$$

for all $\xi \in \mathbb{R}^n \setminus \mathcal{D}$, all $\mu \in \Omega$, where

$$\alpha_3(s,r) = \frac{\rho'(\rho^{-1}(s))\widehat{\alpha}_3(\rho^{-1}(V(\xi))/2, r)}{2}.$$

Since $V$ has local (actually, global) minima at all points in $\mathcal{D}$, it follows that $DV(\xi) \equiv 0$ on $\mathcal{D}$, so we know that the estimate (63) holds also on all of $\mathbb{R}^n$.

Finally, observe that if there exists $\sigma \in \mathcal{K}$ such that (23) holds for all trajectories of the system, then (45) holds for all $\xi$, which, in turn, implies that

$$c_1|h(\xi)| \leq W(\xi) \leq c_2\sigma(|h(\xi)|), \qquad \forall \xi \in \mathbb{R}^n. \tag{64}$$

This results in the desired inequality:

$$\alpha_1(|h(\xi)|) \leq V(\xi) \leq \sigma_1(|h(\xi)|), \qquad \forall \xi \in \mathbb{R}^n, \tag{65}$$

where $\sigma_1(s) = \rho(2c_2\sigma(s))$. This shows that if (23) holds for some $\sigma \in \mathcal{K}$, then property (25) can be strengthened to property (27).

Finally, suppose that, in the above proof, one strengthens (22) to (24). Associated to the function $\beta$ there are, as before, functions $\{T_r\}$. Since also we have an estimate as in (22), there are functions $\{T_r\}$ associated to a $\beta$ as in (22); without loss of generality, we will assume that the same $T_r$'s work for both. Thus, we know that, provided $t \geq T_r(s)$, $|y(t, \xi, u)| \leq s$ whenever $|h(\xi)| \leq r$ or $|\xi| \leq r$. The claim stated after Equation(57) holds now for all $r > |h(\xi)|$ (instead of merely if $r > |\xi|$), because (46) can be strenghtened to

$$\omega(x(t_1, \xi, u)) \leq \beta(|h(\xi)|, t_1).$$

We now repeat the above proof to get a function $W(\xi)$ satisfying (64), and corresponding to (60), one has now also

$$DW(\xi)f(\xi, v) \leq -\widehat{\alpha}_3(W(\xi), |h(\xi)|) \leq -\widehat{\alpha}_3\left(W(\xi), \frac{W(\xi)}{c_1}\right),$$

for all $\xi \in \mathbb{R}^n \setminus \mathcal{D}$ and all $v \in \Omega$. Therefore, on $\mathbb{R}^n \setminus \mathcal{D}$,

$$DW(\xi)f(\xi, v) \leq -\alpha_4(W(\xi)),$$

where $\widetilde{\alpha}_4(s) = \widetilde{\alpha}_3(s, s/c_1)$ is a continuous positive definite function. Using the same smoothing argument as earlier, we can show that there is a smooth function $V$ such that (65) holds for some $\sigma_1, \sigma_2 \in \mathcal{K}_\infty$, and (63) can be strengthened to

$$DV(\xi)f(\xi, v) \leq -\widehat{\alpha}_4(V(\xi)) \tag{66}$$

for all $\xi \in \mathbb{R}^n$, all $v \in \Omega$, where $\widehat{\alpha}_4(\cdot)$ is some continuous positive definite function.

Now we modify the function $V$ to get $V_1$ so that $V_1$ satisfies inequalities of type (64) and (66) with $\widehat{\alpha}_4$ replaced by a $\mathcal{K}_\infty$ function $\alpha_5$. For this purpose, let $\rho_0(\cdot)$ be a smooth $\mathcal{K}_\infty$-function such that $\rho_0(s)\widehat{\alpha}_4(s) \geq 1$ for $s \geq 1$, and let

$$\rho_1(s) = e^{\int_0^s \rho_0(s_1)\, ds_1} - 1.$$

Define $V_1(\xi) = \rho_1(V(\xi))$. It holds that

$$\widehat{\alpha}_1(|h(\xi)|) \leq V_1(\xi) \leq \widehat{\alpha}_2(|h(\xi)|), \qquad \forall \xi \in \mathbb{R}^n,$$

where $\widehat{\alpha}_1(s) = \rho_1(\alpha_1(s))$, and $\widehat{\alpha}_2(s) = \rho_1(\alpha_2(s))$, and

$$DV_1(\xi)f(\xi, v) = -(V_1(\xi) + 1)\rho_0(V(\xi))\widehat{\alpha}_4(V(\xi)) \leq -\alpha_5(V_1(\xi)),$$

for all $\xi \in \mathbb{R}^n$, all $v \in \Omega$, where $\alpha_5$ is any $\mathcal{K}_\infty$ function with the property that

$$\alpha_5(\rho_1(s)) \leq (\rho_1(s) + 1)\rho_0(s)\widehat{\alpha}_4(s)$$

for all $s \geq 0$ (such a $\mathcal{K}_\infty$-function exists because $(s+1)\rho_0(s)\widehat{\alpha}_4(s) \geq s$ for all $s \geq 1$). Using $V_1$ as a Lyapunov function, this completes the proof. ∎

# A  Some Facts Regarding $\mathcal{KL}$ Functions

The following simple observation is proved in [27] and will be needed here too.

**Lemma A.1** For any $\mathcal{KL}$-function $\beta$, there exists a family of mappings $\{T_r\}_{r\geq 0}$ such that

- for each fixed $r > 0$, $T_r : \mathbb{R}_{>0} \xrightarrow{\text{onto}} \mathbb{R}_{>0}$ is continuous and strictly decreasing, and $T_0(s) \equiv 0$;

- for each fixed $s > 0$, $T_r(s)$ is strictly increasing as $r$ increases, and is such that $\beta(r, T_r(s)) < s$, and consequently, $\beta(r,t) < s$ for all $t \geq T_r(s)$.

**Lemma A.2** For any $\mathcal{KL}$ function $\beta$, there exist two $\mathcal{K}$ functions $\kappa_1$ and $\kappa_2$ so that

$$\beta(s,t) \geq \frac{\kappa_1(s)}{1 + \kappa_2(t)} \tag{67}$$

for all $s \geq 0$ and all $t \geq 0$.

*Proof.* We assume that $b := \sup_s \beta(s,0) < \infty$ (otherwise, we first find a $\beta_0 \leq \beta$ with that property, and prove the result for $\beta_0$). We define, for all $s \geq 0$ and $t \geq 0$,

$$\widetilde{\beta}(s,t) := \int_t^{t+1} \beta(s,\tau)\, d\tau \,.$$

Note that $\widetilde{\beta}$ is again of class $\mathcal{KL}$, and $\widetilde{\beta}(s,t) \leq \beta(s,t)$ for all $s,t$. Let

$$\alpha(t) := \sup_{s \geq 0} \widetilde{\beta}(s,t)\,.$$

This is finite everywhere, since it is bounded by $b$. Moreover, it is a continuous function, because

$$\widetilde{\alpha}(t) := \int_t^{t+1} \alpha(\tau)\, d\tau$$

where $\alpha$ is the decreasing function (not necessarily strictly) defined by $\alpha(t) := \sup_{s \geq 0} \beta(s,t)$. We will write from now on "$\widetilde{\beta}(\infty, t)$" instead of "$\widetilde{\alpha}(t)$". Finally, we let

$$\rho(x) := \max\{x, 0\}$$

for all $x \in \mathbb{R}$, and introduce the following function:

$$c : \mathbb{R}^2 \to \mathbb{R} : (x,y) \mapsto -\ln \widetilde{\beta}\left(\frac{1}{\rho(x)}, \rho(y)\right) - \rho(-x) - \rho(-y)$$

where we understand $\widetilde{\beta}(\frac{1}{0}, t)$ as $\widetilde{\alpha}(t)$. As in [1], we let $\mathcal{N}$ denote the class of all functions $k : \mathbb{R} \to \mathbb{R}$ that are nondecreasing, continuous, and unbounded below. Note that $c$ is of class $\mathcal{N}$ on each variable separately. (Continuity follows from the continuity of each of $\widetilde{\beta}(\infty, \cdot)$, $\widetilde{\beta}(s, \cdot)$ for each $s \geq 0$, and $\widetilde{\beta}(\cdot, t)$ for each $t \geq 0$ as well as continuity of $\rho$. The nondecreasing property is clear, using that $\widetilde{\beta}(\cdot, t)$ for each $t \geq 0$ and $\rho$ are nondecreasing, and that $\widetilde{\beta}(\infty, \cdot)$ and $\widetilde{\beta}(s, \cdot)$ for each $s \geq 0$ are nonincreasing. Unbounded below follows from the fact that for $x \to -\infty$ we have $c(x, y_0) = a + x$, where $a = \widetilde{\beta}(\infty, \rho(y_0)) - \rho(-y_0)$ and for $y \to -\infty$ we have $c(x_0, y) = a + y$, where $a = -\ln \widetilde{\beta}\left(\frac{1}{\rho(x_0)}, 0\right) - \rho(-x_0)$.

By Proposition 3.4 in [1], there is some $k \in \mathcal{N}$ such that $c(x,y) \leq k(x) + k(y)$ for all $x, y$. So, we can write, after using that $\beta \geq \widetilde{\beta}$: $\beta(1/x, y) \geq e^{-k(x)}e^{-k(y)}$ for all $x, y > 0$. Equivalently,

$$\beta(s,t) \geq \frac{\kappa_1(s)}{1+\kappa_2(t)}$$

for all $s, t > 0$, when we define

$$\kappa_1(s) := e^{-k(1/s)-k(0)}$$

for all $s > 0$ and

$$\kappa_2(t) := e^{k(t)-k(0)} - 1$$

for all $t \geq 0$. Observe that both of these functions are continuous, nonincreasing, and nonnegative. Moreover, $\kappa_2(0) = 0$, so $\kappa_2$ is in $\mathcal{K}$. From the inequality

$$[1 + \kappa_2(0)]\beta(s,0) \geq \kappa_1(s)$$

for all $s > 0$, and the fact that $\beta(0,0) = 0$, we conclude that $\lim_{s \to 0^+} \kappa_1(s) = 0$, so we may extend $\kappa_1$ by defining $\kappa_1(0) = 0$, and thus $\kappa_1$ is in $\mathcal{K}$ as well. ∎

As $\kappa_1$ and $\kappa_2$ in Lemma A.2 are continuous, we have, in particular:

**Corollary A.3** *For any $\mathcal{KL}$-function $\beta$, there is a (jointly) continuous $\mathcal{KL}$-function $\beta_1$ such that $\beta(s,r) \geq \beta_1(s,r)$ for all $(s,r) \in \mathbb{R}_{\geq 0} \times \mathbb{R}_{\geq 0}$.* □

The following is a generalization of the comparison lemma given in [13]. It plays a role in the proofs of sufficiency, which are the easier parts of the theorems.

**Lemma A.4** *For any $\mathcal{K}$-function $\kappa$, there exists a $\mathcal{KL}$ function $\beta$ such that if $y(\cdot)$ is any locally absolutely continuous function defined on some interval $[0, T]$ with $y(t) \geq 0$, and if $y(\cdot)$ satisfies the differential inequality*

$$\dot{y}(t) \leq -c\kappa(y(t)) \quad \text{for almost all} \quad t \in [0,T] \tag{68}$$

*for some $c \geq 0$ with $y(0) = y_0 \geq 0$, then it holds that*

$$y(t) \leq \beta(y_0, ct)$$

*for all $t \in [0, T]$.*

*Proof.* First of all, by Lemma 4.4 in [13], for each $\kappa \in \mathcal{K}$, there exists $\beta \in \mathcal{KL}$ such that for any locally absolutely continuous function $z(t) \geq 0$, if it satisfies the inequality

$$\dot{z}(t) \leq -\kappa(z(t))$$

on $[0, T]$, it holds that $z(t) \leq \beta(z(0), t)$ for all $t$. (The statement in that reference applies to $z$ defined on all of $[0, \infty)$, but exactly the same proof works for a finite interval.)

Let $y(t)$ be a function as in the statement of the lemma for some $c > 0, T > 0$. Let $\widetilde{y}(t) = y(t/c)$. Then $\widetilde{y}$ is again locally absolutely continuous and nonnegative on $[0, cT]$. Moreover, $\widetilde{y}$ satisfies the inequality

$$\frac{d}{dt}\widetilde{y}(t) \leq -\kappa(\widetilde{y}(t)).$$

Hence,
$$\widetilde{y}(t) \leq \beta(\widetilde{y}(0), t)$$
for all $t \in [0, cT]$. This then implies that
$$y(t) \leq \beta(y(0), ct)$$
for all $t \in [0, T]$. ∎

Finally, we have this fact, mentioned when discussing decrease conditions:

**Lemma A.5** Let $V : \mathbb{R}^n \to \mathbb{R}$ be a $C^1$ positive definition function with the following property: for some $\mathcal{K}$ function $\chi$, it holds that
$$V(\xi) \geq \chi(|\mu|) \text{ and } V(\xi) \neq 0 \;\Rightarrow\; DV(\xi)f(\xi,\mu) < 0.$$
Then, there is a function $\alpha \in \mathcal{KL}$ so that:
$$V(\xi) \geq \chi(|\mu|) \;\Rightarrow\; DV(\xi)f(\xi,\mu) \leq -\alpha(V(\xi), |\xi|)$$
for all $\xi \in \mathbb{R}^n, \mu \in \mathbb{R}^m$.

*Proof.* Without loss of generality, we assume that $\chi \in \mathcal{K}_\infty$. Define the set, for each $s, t \geq 0$:
$$R(s,t) := \{(x,u) : |\xi| \leq t, \; V(\xi) \geq s, \; |\mu| \leq \chi^{-1}(V(\xi))\}.$$
These sets are compact (possibly empty), for each $s$ and $t$. Note the following properties:
$$s > s' \;\Rightarrow\; R(s,t) \subseteq R(s',t)$$
$$t > t' \;\Rightarrow\; R(s,t') \subseteq R(s,t).$$
Now let
$$\alpha_0(s,t) = \min_{(\xi,\mu)\in R(s,t)} -DV(\xi)f(\xi,\mu)$$
(with the convention that $\alpha_0(s,t) = +\infty$ if $R(s,t) = \emptyset$). Then, $\alpha_0(s,t)$ is nonincreasing in $t$ and nondecreasing in $s$. Moreover, $\alpha(s,t) > 0$ whenever $s > 0$ (by the hypothesis of the lemma). Next let:
$$\widehat{\alpha}(s,t) := \min\{\alpha_0(s,t), s\}.$$
This function has the same monotonicity properties as $\alpha_0$, it satisfies $\alpha_0(s,t) \geq \widehat{\alpha}(s,t)$ for all $s,t$, and is finite-valued. It also satisfies $\widehat{\alpha}(s,t) \neq 0$ for $s > 0$. Now pick
$$\widetilde{\alpha}(s,t) := \int_{s-1}^{s} \widehat{\alpha}(\sigma, t)\,d\sigma$$
(let $\widehat{\alpha}(s,t) := 0$ for $s < 0$). This function still has the same monotonicity properties, satisfies $\widetilde{\alpha}(s,t) > 0$ for $s > 0$, and it is continuous in $s$. It may not be strictly increasing in $s$, nor need it converge to zero as $t \to 0$, so we obtain finally a $\mathcal{KL}$ function $\alpha$ by defining:
$$\alpha(s,t) := \frac{s\widetilde{\alpha}(s,t)}{(1+s)(1+t)}.$$
This satisfies the desired properties by construction, because
$$V(\xi) \geq \chi(|\mu|) \;\Rightarrow\; DV(\xi)f(\xi,\mu) \leq -\alpha(V(\xi), |\mu|)$$
and $\alpha_0 \geq \widehat{\alpha} \geq \widetilde{\alpha} \geq \alpha$ pointwise. ∎